\documentclass[a4paper,oneside,notitlepage]{article}%
\usepackage{amsmath}
\usepackage{amssymb}
\usepackage{amscd}
\usepackage{amsfonts}
\usepackage{amstext}
\usepackage{graphicx}%
\newtheorem{theorem}{Theorem}

\newtheorem{remark}[theorem]{Remark}

\newtheorem{definition}[theorem]{Definition}

\ifx\pdfoutput\relax\let\pdfoutput=\undefined\fi
\newcount\msipdfoutput
\ifx\pdfoutput\undefined\else
\ifcase\pdfoutput\else
\ifx\paperwidth\undefined\else
\ifdim\paperheight=0pt\relax\else\pdfpageheight\paperheight\fi
\ifdim\paperwidth=0pt\relax\else\pdfpagewidth\paperwidth\fi
\fi\fi\fi
\begin{document}

\title{Domain of existence of the Laplace transform of infinitely divisible negative
multinomial distributions}
\author{Philippe Bernardoff\\Universit\'{e} de Pau et des Pays de l'Adour\\Laboratoire de Math\'{e}matiques et de leurs Applications de Pau,\\UMR 5142, BP1155, 64013 PAU France\\e.mail : philippe.bernardoff@univ-pau.fr}
\date{Mai 30, 2021 }
\maketitle

\begin{abstract}
This article provides the domain of existence $\Omega$ of the Laplace
transform of infinitely divisible negative multinomial distributions. We
define a negative multinomial distribution on $\mathbb{N}^{n},$ where
$\mathbb{N}$ is the set of nonnegative integers, by its probability generating
function which will be of the form \newline$\left(  A\left(  a_{1}z_{1}%
,\ldots,a_{n}z_{n}\right)  /A\left(  a_{1},\ldots,a_{n}\right)  \right)
^{-\lambda}$ where $A\left(  \mathbf{z}\right)  =%
%TCIMACRO{\dsum \limits_{T\subset\left\{  1,2,\ldots,n\right\}  }}%
%BeginExpansion
{\displaystyle\sum\limits_{T\subset\left\{  1,2,\ldots,n\right\}  }}
%EndExpansion
a_{T}\prod\limits_{i\in T}z_{i},$ where $a_{\emptyset}\neq0,$ and where
$\lambda$ is a positive number. Finding couples $\left(  A,\lambda\right)  $
for which we obtain a probability generating function is a difficult problem.
Necessary and sufficient conditions on the coefficients $a_{T}$ of $A$ for
which we obtain a probability generating function for any positive number
$\lambda$ are know by (Bernardoff, 2003). Thus we obtain necessary and
sufficient conditions on $\mathbf{a}=\left(  a_{1},\ldots,a_{n}\right)  $ so
that $\mathbf{a}=\left(  \mathbf{e}^{t_{1}},\ldots,\mathbf{e}^{t_{n}}\right)
$ with $\mathbf{t}=\left(  t_{1},\ldots,t_{n}\right)  $ belonging to $\Omega$.
This makes it possible to construct all the infinitely divisible multinomial
distributions on $\mathbb{N}^{n}$. We give examples of construction in
dimensions 2 and 3.


\newpage

\end{abstract}

\section{Introduction}

In this article, we consider the following definition, see references in
Bernardoff (2003). We shall say that the probability distribution
$\sum\limits_{\mathbf{\boldsymbol{\alpha}}\in\mathbb{N}^{n}}%
p_{\mathbf{\boldsymbol{\alpha}}}\delta_{\mathbf{\boldsymbol{\alpha}}}$ on
$\mathbb{N}^{n}$, where $n$ is a non negative number, is a \textit{negative
multinomial distribution} if there exists an affine polynomial $P\left(
z_{1},\ldots,z_{n}\right)  $ and $\lambda>0$ such that $P\left(
0,\ldots,0\right)  \neq0$, and $P\left(  1,\ldots,1\right)  =1.$%
\begin{equation}
\sum\limits_{\boldsymbol{\alpha}\in\mathbb{N}^{n}}p_{\boldsymbol{\alpha}}%
z_{1}^{\alpha_{1}}\cdots z_{n}^{\alpha_{n}}=\left(  P\left(  z_{1}%
,\ldots,z_{n}\right)  \right)  ^{-\lambda}. \label{zalpha}%
\end{equation}
That means a polynomial which is affine with respect to each $z_{j},$
$j=1,\cdots,n,$ or for which $\dfrac{\partial^{2}}{\partial z_{j}^{2}}P=0$ for
all $j=1,\ldots,n.$ That is $P\left(  z_{1},\ldots,z_{n}\right)
=\sum\limits_{T\in\mathfrak{P}_{n}}a_{T}\mathbf{z}^{T},$ where $\mathfrak{P}%
_{n}$ is the set of the subset of $\left\{  1,2,\ldots,n\right\}  =\left[
n\right]  ,$ and where $\mathbf{z}^{T}=\prod_{t\in T}z_{t}.$if $\mathbf{z}%
=\left(  z_{1},\ldots,z_{n}\right)  \in\mathbb{R}^{n}.$ For instance, for
$n=2,$ such as $P$ has the form $P\left(  z_{1},z_{2}\right)  =a_{\emptyset
}+a_{\left\{  1\right\}  }z_{1}+a_{\left\{  2\right\}  }z_{2}+a_{\left\{
1,2\right\}  }z_{1}z_{2}$ with $a_{\emptyset}\neq0.$ However, finding exactly
which pairs $\left(  P,\lambda\right)  $ are compatible is an unsolved
problem.\newline Before giving the main result, let us make an observation. If
\textbf{$\boldsymbol{\alpha}$}$=\left(  \alpha_{1},\ldots,\alpha_{n}\right)
\in\mathbb{N}^{n},$ then we denote%
\begin{equation}
\mathbf{z}^{\mathbf{\boldsymbol{\alpha}}}=\prod_{i=1}^{n}z_{i}^{\alpha_{i}%
}=z_{1}^{\alpha_{1}}\ldots z_{n}^{\alpha_{n}}. \label{zalpha2}%
\end{equation}
Let $A$ be any polynomial such that $A\left(  0,\ldots,0\right)  =1,$ and
suppose that the Taylor expansion%
\[
\left(  A\left(  z_{1},\ldots,z_{n}\right)  \right)  ^{-\lambda}%
=\sum_{\mathbf{\boldsymbol{\alpha}}\in\mathbb{N}^{n}}%
c_{\mathbf{\boldsymbol{\alpha}}}\left(  \lambda\right)  \mathbf{z}%
^{\mathbf{\boldsymbol{\alpha}}}%
\]
has non-negative coefficients $c_{\mathbf{\boldsymbol{\alpha}}}\left(
\lambda\right)  .$ Let $a_{1,}\ldots,a_{n}$ be positive numbers such that
$\sum_{\mathbf{\boldsymbol{\alpha}}\in\mathbb{N}^{n}}%
c_{\mathbf{\boldsymbol{\alpha}}}\left(  \lambda\right)  a_{1}^{\alpha_{1}%
}\ldots a_{n}^{\alpha_{n}}<\infty.$ With such a sequence $\mathbf{a}=\left(
a_{1},\ldots,a_{n}\right)  $ we associate the negative multinomial
distribution $\sum\limits_{\mathbf{\boldsymbol{\alpha}}\in\mathbb{N}^{n}%
}p_{\mathbf{\boldsymbol{\alpha}}}\delta_{\mathbf{\boldsymbol{\alpha}}}$
defined by
\begin{equation}
\sum\limits_{\mathbf{\boldsymbol{\alpha}}\in\mathbb{N}^{n}}%
p_{\mathbf{\boldsymbol{\alpha}}}\mathbf{z}^{\mathbf{\boldsymbol{\alpha}}%
}=\left(  \dfrac{A\left(  a_{1}z_{1},\ldots,a_{n}z_{n}\right)  }{A\left(
a_{1},\ldots,a_{n}\right)  }\right)  ^{-\lambda}, \label{P-lambda}%
\end{equation}
thus%
\begin{equation}
P\left(  z_{1},\ldots,z_{n}\right)  =\frac{A\left(  a_{1}z_{1},\ldots
,a_{n}z_{n}\right)  }{A\left(  a_{1},\ldots,a_{n}\right)  } \label{Pofzalpha}%
\end{equation}
in the notation (\ref{zalpha}).

Bernardoff (2003) define the polynomials $b_{T}$ by

\begin{definition}
Let $P\left(  \mathbf{z}\right)  =\sum\limits_{T\in\mathcal{P}_{n}}%
a_{T}\mathbf{z}^{T}$ be an affine polynomial $P\left(  z_{1},\ldots
,z_{n}\right)  $ such that $P\left(  0,\ldots,0\right)  =0,$ and $A=1-P$. Let
$T$ be in $\mathfrak{P}_{n}^{\ast}$ the set of the nonempty subset of $\left[
n\right]  $ let us denote by $b_{T}$ the number defined by
\[
b_{T}=\left.  \frac{\partial^{\left\vert T\right\vert }}{\partial z^{T}%
}\left(  \log\left(  1-P\right)  \right)  \right\vert _{\mathbf{0}},
\]
where then $\left\vert T\right\vert $ is the cardinal of $T$ and $\partial
z^{T}=\prod_{t\in T}\partial z_{t},$ then
\begin{equation}
b_{T}=\sum\limits_{l=1}^{\left\vert T\right\vert }\left(  l-1\right)
!\sum\limits_{\mathcal{T}\in\Pi_{T}^{l}}a_{\mathcal{T}} \label{bT}%
\end{equation}
where $\Pi_{T}$ is the set of the partition of $T$, and $\Pi_{T}^{l}$ is the
set of the partition of lenght $l$ of $T$ (if $\mathcal{T}=\left\{
T_{1},T_{2},\ldots,T_{l}\right\}  ,$ the partition $\mathcal{T}$ of $T$ is of
length $l$).
\end{definition}

For instance, for $n=3,b_{\left\{  1\right\}  }=a_{\left\{  1\right\}  },$
$b_{\left\{  1,2\right\}  }=a_{\left\{  1,2\right\}  }+a_{\left\{  1\right\}
}a_{\left\{  2\right\}  }$ and $b_{\left\{  1,2,3\right\}  }=a_{\left\{
1,2,3\right\}  }+a_{\left\{  1\right\}  }a_{\left\{  2,3\right\}
}+a_{\left\{  2\right\}  }a_{\left\{  1,3\right\}  }+a_{\left\{  3\right\}
}a_{\left\{  1,2\right\}  }+2a_{\left\{  1\right\}  }a_{\left\{  2\right\}
}a_{\left\{  3\right\}  }$. Now, if there is no ambiguity, for simplicity we
omit the braces.

Using the numbers $b_{T}$, Bernardoff (2003) proves the following theorem.

\begin{theorem}
\label{theoremID}Let $P\left(  z\right)  =%
%TCIMACRO{\dsum \nolimits_{T\in\mathfrak{P}_{n}\ast}}%
%BeginExpansion
{\displaystyle\sum\nolimits_{T\in\mathfrak{P}_{n}\ast}}
%EndExpansion
a_{T}z^{T},$ as before, and suppose that $\left(  1-P\left(  z\right)
\right)  ^{-\lambda}=%
%TCIMACRO{\dsum \nolimits_{\mathbf{\boldsymbol{\alpha}}\in\mathbb{N}^{n}}}%
%BeginExpansion
{\displaystyle\sum\nolimits_{\mathbf{\boldsymbol{\alpha}}\in\mathbb{N}^{n}}}
%EndExpansion
c_{\mathbf{\boldsymbol{\alpha}}}\left(  \lambda\right)
z^{\mathbf{\boldsymbol{\alpha}}}$. Then $c_{\mathbf{\boldsymbol{\alpha}}%
}\left(  \lambda\right)  \geqslant0\ $for all positive $\lambda$ if and only
if $b_{T},$ given by $\left(  \ref{bT}\right)  $, is non negative for all
$T\in\mathfrak{P}_{n}^{\ast}.$
\end{theorem}

See examples in dimension $n=2,3$ in Bernardoff (2003).

This article is organized as follows. Section 2 gives the main result. Section
3 applies the main result to bivariate and trivariate cases.

\section{Domain of existence of the Laplace transform}

Let $A$ be an affine polynomial on $\mathbb{R}^{n}$ and let $\lambda>0$ such
that $A\left(  0,\ldots,0\right)  =1$ and such that $\left(  A\left(
z_{1},\ldots,z_{n}\right)  \right)  ^{-\lambda}=%
%TCIMACRO{\dsum \limits_{\mathbf{\boldsymbol{\alpha}}\in\mathbb{N}^{n}}}%
%BeginExpansion
{\displaystyle\sum\limits_{\mathbf{\boldsymbol{\alpha}}\in\mathbb{N}^{n}}}
%EndExpansion
c_{\mathbf{\boldsymbol{\alpha}}}\left(  \lambda\right)
z^{\mathbf{\boldsymbol{\alpha}}}$ satisfies $c_{\mathbf{\boldsymbol{\alpha}}%
}\left(  \lambda\right)  >0$ for all \textbf{$\boldsymbol{\alpha}$} in
$\mathbb{N}^{n}.$ Consider the discrete measure on $\mathbb{N}^{n}$,
$\mu_{\lambda}=%
%TCIMACRO{\dsum \limits_{\mathbf{\boldsymbol{\alpha}}\in\mathbb{N}^{n}}}%
%BeginExpansion
{\displaystyle\sum\limits_{\mathbf{\boldsymbol{\alpha}}\in\mathbb{N}^{n}}}
%EndExpansion
c_{\mathbf{\boldsymbol{\alpha}}}\left(  \lambda\right)  \delta
_{\mathbf{\boldsymbol{\alpha}}}$.

The present section aims to describe the convex set
\[
D\left(  \mu_{\lambda}\right)  =\left\{  \boldsymbol{\theta}\in\mathbb{R}^{n},%
%TCIMACRO{\dsum \limits_{\mathbf{\boldsymbol{\alpha}}\in\mathbb{N}^{n}}}%
%BeginExpansion
{\displaystyle\sum\limits_{\mathbf{\boldsymbol{\alpha}}\in\mathbb{N}^{n}}}
%EndExpansion
c_{\mathbf{\boldsymbol{\alpha}}}\left(  \lambda\right)  \mathbf{e}^{\theta
_{1}\alpha_{1}+\cdots+\theta_{n}\alpha_{n}}<+\infty\right\}  ,
\]
which is an important object in order to study the natural exponential family
generated by $\mu_{\lambda}$ (see Letac, 1991 and Bar-Lev \textit{et al.},
1994). The answer is contained in the following Proposition.

\begin{theorem}
\label{mainTheorem}With the above notation, denote $H=\left\{  \mathbf{s}%
\in\mathbb{R}^{n},s_{1}+\cdots+s_{n}=0\right\}  .$ For \textbf{$s$}$\in H,$ we
denote by $R_{\mathbf{s}}$ the smallest positive zero of the polynomial
$P_{\mathbf{s}}\left(  t\right)  =A\left(  t\mathbf{e}^{s_{1}},\ldots
,t\mathbf{e}^{s_{n}}\right)  $. Then the map $\mathbf{s}\mapsto\mathbf{s}+\log
R_{\mathbf{s}}\left(  1,\ldots,1\right)  $ is a parametrization by $H$ of a
hypersurface in $\mathbb{R}^{n}$ which is the boundary of $D\left(
\mu_{\lambda}\right)  .$ More specifically if $\boldsymbol{\theta}$ is in
$\mathbb{R}^{n}$ and if $\overline{\theta}_{n}=\dfrac{1}{n}\left(  \theta
_{1}+\cdots+\theta_{n}\right)  $ and if $\mathbf{s}=\boldsymbol{\theta
}-\overline{\theta}_{n}\left(  1,\ldots,1\right)  ,$ then $\boldsymbol{\theta
}$ is in $D\left(  \mu_{\lambda}\right)  $ if and only if $\overline{\theta
}_{n}<\log R_{\mathbf{s}}.$\newline Finally $D\left(  \mu_{\lambda}\right)  $
is an open set.
\end{theorem}

\begin{description}
\item[Proof] We first prove that the radius of convergence $R$ of the power
series%
\begin{equation}
P_{\mathbf{s}}^{-\lambda}\left(  t\right)  =%
%TCIMACRO{\dsum \limits_{n\in\mathbb{N}}}%
%BeginExpansion
{\displaystyle\sum\limits_{n\in\mathbb{N}}}
%EndExpansion
u_{n}\left(  \lambda\right)  t^{n} \label{radiusR}%
\end{equation}
\newline is equal to $R_{\mathbf{s}}.$ This comes from the following fact:
since $u_{n}\left(  \lambda\right)  \geqslant0,$ a known result in theory of
analytic functions (see Titchmarsh (1939) \textbf{7.21}) implies that
$t\mapsto P_{\mathbf{s}}^{-\lambda}\left(  t\right)  $ is not analytic at $R.$
Since $P_{\mathbf{s}}\left(  0\right)  =1,$ $P_{\mathbf{s}}\left(  t\right)
>0$ for $0<t<R_{\mathbf{s}}$ and $P_{\mathbf{s}}\left(  R_{\mathbf{s}}\right)
=0$ clearly $R=R_{\mathbf{s}}.$

We now observe that if \textbf{$\boldsymbol{\theta}$}$=\left(  \theta
_{1},\ldots,\theta_{n}\right)  $ is such that
\[%
%TCIMACRO{\dsum \limits_{\mathbf{\boldsymbol{\alpha}}\in\mathbb{N}^{n}}}%
%BeginExpansion
{\displaystyle\sum\limits_{\mathbf{\boldsymbol{\alpha}}\in\mathbb{N}^{n}}}
%EndExpansion
c_{\mathbf{\boldsymbol{\alpha}}}\left(  \lambda\right)  \mathbf{e}^{\alpha
_{1}\theta_{1}+\cdots+\alpha_{n}\theta_{n}}<+\infty
\]
\newline then for all $p\geqslant0$ we have
\[%
%TCIMACRO{\dsum \limits_{\mathbf{\boldsymbol{\alpha}}\in\mathbb{N}^{n}}}%
%BeginExpansion
{\displaystyle\sum\limits_{\mathbf{\boldsymbol{\alpha}}\in\mathbb{N}^{n}}}
%EndExpansion
c_{\mathbf{\boldsymbol{\alpha}}}\left(  \lambda\right)  \mathbf{e}^{\alpha
_{1}\left(  \theta_{1}-p\right)  +\cdots+\alpha_{n}\left(  \theta
_{n}-p\right)  }<+\infty.
\]
\newline We now fix $\boldsymbol{\theta}$ in $D\left(  \mu_{\lambda}\right)
.$ Write $\overline{\theta}_{n}=\left(  \theta_{1}+\cdots+\theta_{n}\right)
/n.$ The orthogonal projection of $\boldsymbol{\theta}$ on $H$ is
$\mathbf{s}=\boldsymbol{\theta}-\overline{\theta}_{n}\left(  1,\ldots
,1\right)  =\left(  s_{1},\ldots,s_{n}\right)  $. Thus for all $j=1,\ldots,n$
we have $\theta_{j}-s_{j}=\overline{\theta}_{n}.$ We claim that $\overline
{\theta}_{n}<\log R_{\mathbf{s}}.$ If not, we have $t_{0}=$\textbf{$e$%
}$^{\overline{\theta}_{n}}\geqslant R_{\mathbf{s}}$. But $A\left(
\mathbf{e}^{\theta_{1}},\ldots\mathbf{e}^{\theta_{n}}\right)  $ is
$P_{\mathbf{s}}\left(  t_{0}\right)  $ and the previous remark shows that for
all $p\geqslant0,$ $p\mapsto P_{\mathbf{s}}\left(  \mathbf{e}^{-p}%
t_{0}\right)  $ is always positive. This contradicts the fact that
$t_{0}\geqslant R_{\mathbf{s}}$.\newline Conversely if $\boldsymbol{\theta}$
is such that $\overline{\theta}_{n}<\log R_{\mathbf{s}}$ with $\mathbf{s}%
=\boldsymbol{\theta}-\overline{\theta}_{n}\left(  1,\ldots,1\right)  $ a
similar reasoning shows that $\boldsymbol{\theta}\in D\left(  \mu_{\lambda
}\right)  .$\newline Finally for $t=R_{\mathbf{s}}$ in $\left(  \ref{radiusR}%
\right)  $ the series diverges. A short proof goes as follows:\newline Write
$P_{\mathbf{s}}\left(  t\right)  =\left(  1-\dfrac{t}{r_{0}}\right)  \left(
1-\dfrac{t}{r_{1}}\right)  \cdots\left(  1-\dfrac{t}{r_{k}}\right)  $ where
$\left\vert r_{j}\right\vert \geqslant r_{0}=R_{\mathbf{s}}$ by definition of
$R_{\mathbf{s}}.$\newline Thus $%
%TCIMACRO{\dsum \limits_{n\in\mathbb{N}}}%
%BeginExpansion
{\displaystyle\sum\limits_{n\in\mathbb{N}}}
%EndExpansion
u_{n}\left(  \lambda\right)  t^{n}$ is the product of Newton Series $%
%TCIMACRO{\dsum \limits_{n\in\mathbb{N}}}%
%BeginExpansion
{\displaystyle\sum\limits_{n\in\mathbb{N}}}
%EndExpansion
\dfrac{1}{n!}\left\langle \lambda\right\rangle _{n}\left(  \dfrac{t}{r_{k}%
}\right)  ^{n}$ and the series corresponding to $k=0$ diverges for $t=r_{0}.$
Thus $%
%TCIMACRO{\dsum \limits_{n\in\mathbb{N}}}%
%BeginExpansion
{\displaystyle\sum\limits_{n\in\mathbb{N}}}
%EndExpansion
u_{n}\left(  \lambda\right)  r_{0}^{n}=+\infty.\blacksquare$

\begin{remark}
\label{mainremark}With the notations of Theorem \ref{mainTheorem}, if
$\boldsymbol{\theta}$ is in the boundary of $D\left(  \mu_{\lambda}\right)  $,
$\exists\mathbf{s}\in H,\boldsymbol{\theta}=\mathbf{s}+\log R_{\mathbf{s}%
}\left(  1,\ldots,1\right)  ,$ then $z_{i}=$\textbf{$e$}$^{\theta_{i}%
}=R_{\mathbf{s}}$\textbf{$e$}$^{s_{i}}$ for $i\in\left[  n\right]  ,$ and
\[
A\left(  z_{1},\ldots,z_{n}\right)  =A\left(  R_{\mathbf{s}}\mathbf{e}^{s_{1}%
},\ldots,R_{\mathbf{s}}\mathbf{e}^{s_{n}}\right)  =0,
\]
by the definition of $R_{\mathbf{s}}.$
\end{remark}
\end{description}

\section{Examples in dimensions 2 and 3}

\textbf{Example 1. }For $n=2,$ we take $a_{1}=1,$ $a_{2}=1$ and $a_{1,2}%
=a\geqslant-1$ so that the conditions of Theorem \ref{theoremID} are
satisfied. Hence $A\left(  z_{1},z_{2}\right)  =1-z_{1}-z_{2}-az_{1}z_{2}.$
For $z_{1}=t\mathbf{e}^{s_{1}},$ $z_{2}=t\mathbf{e}^{s_{2}},$ with
$s_{1}+s_{2}=0$ we have $A\left(  z_{1},z_{2}\right)  =1-t\left(
\mathbf{e}^{s_{1}}+\mathbf{e}^{-s_{1}}\right)  -at^{2}.$ If \textbf{$s$%
}$=\left(  s_{1},-s_{1}\right)  ,$ $R_{\mathbf{s}}=\frac{1}{a}\left(  -\cosh
s_{1}+\sqrt{\frac{1}{2}+a+\frac{1}{2}\cosh2s_{1}}\right)  ,$ we obtain the
parametrization of the boundary of $D\left(  u_{\lambda}\right)  $ :
\[
\left\{
\begin{tabular}
[c]{l}%
$x=\theta_{1}=s_{1}+\log\frac{1}{a}\left(  -\cosh s_{1}+\sqrt{\frac{1}%
{2}+a+\frac{1}{2}\cosh2s_{1}}\right)  $\\
$y=\theta_{2}=-s_{1}+\log\frac{1}{a}\left(  -\cosh s_{1}+\sqrt{\frac{1}%
{2}+a+\frac{1}{2}\cosh2s_{1}}\right)  $%
\end{tabular}
\ \right.  ,
\]
whose graphic representation is%

\begin{center}
\includegraphics{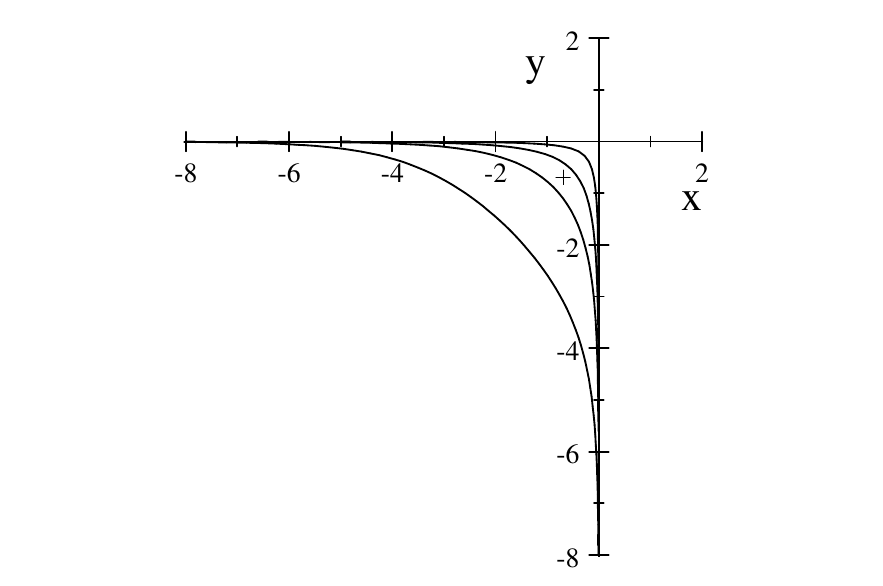}
\\
figure 1 : the boundary of $D\left(  u_{\lambda}\right)  $ for $a=-\frac
{9}{10},-\frac{1}{2},1,20.$%
\label{fig1}%

\end{center}
%EndExpansion

\begin{remark}
Using the Remark \ref{mainremark} we obtain an other parametrization of the
boundary of $D\left(  u_{\lambda}\right)  $ :
\[
\theta_{1}<0,\theta_{2}=-\log\left(  1+\left(  a+1\right)  /\left(
\mathbf{e}^{-\theta_{1}}-1\right)  \right)  .
\]
In addition, because $D\left(  \mu_{\lambda}\right)  $ is a convex set,
$\mathbf{\theta}=\left(  \theta_{1},\theta_{2}\right)  \in D\left(
\mu_{\lambda}\right)  $ is defined by
\begin{equation}
\theta_{1}<0,\theta_{2}<-\log\left(  1+\left(  a+1\right)  /\left(
\mathbf{e}^{-\theta_{1}}-1\right)  \right)  . \label{conditiondim2exemple}%
\end{equation}

\end{remark}

In this case, for $a=-1/2,$ $R_{\mathbf{s}}=2\cosh s_{1}-\sqrt{2\cosh2s_{1}}$
and if we choose $s_{1}=0,$ the condition $\overline{\theta}_{2}<\log
R_{\mathbf{s}}$ becomes $\overline{\theta}_{2}<-\log\left(  1+\sqrt
{2}/2\right)  .$ As $\overline{\theta}_{2}=-\log2<-\log\left(  1+\sqrt
{2}/2\right)  $, then $\left(  -\log2,-\log2\right)  \in D\left(  \mu
_{\lambda}\right)  $. Then, the introduction proves that $(A\left(  \frac
{1}{2}z_{1},\frac{1}{2}z_{2}\right)  /A\left(  \frac{1}{2},\frac{1}{2}\right)
)^{-\lambda}=\left(  8-4z_{1}-4z_{2}+z_{1}z_{2}\right)  ^{-\lambda}$ is a
generating function for all $\lambda>0$.

\begin{remark}
For $a=-1/2,$ the condition \ref{conditiondim2exemple} gives for $\theta
_{1}=-\log2,$ $\theta_{2}<-\log\left(  3/2\right)  ,$ and $\theta_{2}=-\log2$
is suitable. Hence $\left(  -\log2,-\log2\right)  \in D\left(  \mu_{\lambda
}\right)  .$
\end{remark}

Again with $a=-1/2,$ if we choose $s_{1}=\log2,$ $R_{\mathbf{s}}=\frac{5}%
{2}-\frac{1}{2}\sqrt{17},$ then $\overline{\theta}_{2}=-2\log2<\log
R_{\mathbf{s}}$ and $\mathbf{\theta}=\mathbf{s}+\overline{\theta}_{2}\left(
1,1\right)  =\allowbreak\left(  -\log2,-3\log2\right)  $. Then, the
introduction proves that $(A\left(  \frac{1}{2}z_{1},\frac{1}{8}z_{2}\right)
/A\left(  \frac{1}{2},\frac{1}{8}\right)  )^{-\lambda}=\left(  \frac{32}%
{13}-\frac{16}{13}z_{1}-\frac{4}{13}z_{2}+\frac{1}{13}z_{1}z_{2}\right)
^{-\lambda}$ is a generating function for all $\lambda>0$.

\begin{remark}
For $a=-1/2,$ the condition \ref{conditiondim2exemple} gives for $\theta
_{1}=-\log2,$ $\theta_{2}<-\log\left(  3/2\right)  ,$ and $\theta_{2}=-3\log2$
is suitable. Hence $\left(  -\log2,-3\log2\right)  \in D\left(  \mu_{\lambda
}\right)  .$ This method is easier to use.
\end{remark}

\textbf{Example 2. }For $n=3,$ the conditions of Theorem \ref{theoremID} are
for $i,j=1,2,3$:\newline$b_{i}=a_{i}\geqslant0$ ; $a_{ij}\geqslant-a_{i}a_{j}$
; $a_{123}\geqslant-\left(  a_{1}a_{23}+a_{2}a_{13}+a_{3}a_{12}+2a_{1}%
a_{2}a_{3}\right)  $\newline We take $a_{1}=a_{2}=a_{3}=1,$ $a_{12}%
=a_{13}=a_{23}=a$ and $a_{123}=b$, so that
\[
A\left(  z\right)  =1-\left(  \left(  z_{1}+z_{2}+z_{3}\right)  +a\left(
z_{1}z_{2}+z_{1}z_{3}+z_{2}z_{3}\right)  +bz_{1}z_{2}z_{3}\right)  .
\]
The conditions of Theorem \ref{theoremID} are satisfied for $a\geqslant-1$ and
$b\geqslant-3a-2.$ We take $a=1$ and $b=0,$ hence
\[
A\left(  z\right)  =1-z_{1}-z_{2}-z_{3}-z_{1}z_{2}-z_{1}z_{3}-z_{2}z_{3}.
\]
Let $z_{1}=t\mathbf{e}^{s_{1}},$ $z_{2}=t\mathbf{e}^{s_{2}}$ and
$z_{3}=t\mathbf{e}^{s_{2}},$ with $s_{1}+s_{2}+s_{3}=0,$ then $P_{\mathbf{s}%
}\left(  t\right)  =1-\left(  \mathbf{e}^{s_{1}}+\mathbf{e}^{s_{2}}%
+\mathbf{e}^{-s_{1}-s_{2}}\right)  t-\left(  \mathbf{e}^{-s_{1}}%
+\mathbf{e}^{-s_{2}}+\mathbf{e}^{s_{1}+s_{2}}\right)  t^{2}=0,$ and we have
\[
R_{\mathbf{s}}=\frac{-\mathbf{e}^{s_{1}}-\mathbf{e}^{s_{2}}-\mathbf{e}%
^{-s_{1}-s_{2}}+\sqrt{\left(  \mathbf{e}^{2s_{1}}+\mathbf{e}^{2s_{2}%
}+6\mathbf{e}^{s_{1}+s_{2}}+\mathbf{e}^{-2(s_{1}+s_{2})}+6\mathbf{e}^{-s_{1}%
}+6\mathbf{e}^{-s_{2}}\right)  }}{2\left(  \mathbf{e}^{-s_{1}}+\mathbf{e}%
^{-s_{2}}+\mathbf{e}^{s_{1}+s_{2}}\right)  }.
\]
Finally, the parametrization of the boundary of $D\left(  u_{\lambda}\right)
$ is:%
\[
\left\{
\begin{tabular}
[c]{l}%
$x=\theta_{1}=s_{1}+\log\frac{-\mathbf{e}^{s_{1}}-\mathbf{e}^{s_{2}%
}-\mathbf{e}^{-s_{1}-s_{2}}+\sqrt{\left(  \mathbf{e}^{2s_{1}}+\mathbf{e}%
^{2s_{2}}+6\mathbf{e}^{s_{1}+s_{2}}+\mathbf{e}^{-2(s_{1}+s_{2})}%
+6\mathbf{e}^{-s_{1}}+6\mathbf{e}^{-s_{2}}\right)  }}{2\left(  \mathbf{e}%
^{-s_{1}}+\mathbf{e}^{-s_{2}}+\mathbf{e}^{s_{1}+s_{2}}\right)  }$\\
$y=\theta_{2}=s_{2}+\log\frac{-\mathbf{e}^{s_{1}}-\mathbf{e}^{s_{2}%
}-\mathbf{e}^{-s_{1}-s_{2}}+\sqrt{\left(  \mathbf{e}^{2s_{1}}+\mathbf{e}%
^{2s_{2}}+6\mathbf{e}^{s_{1}+s_{2}}+\mathbf{e}^{-2(s_{1}+s_{2})}%
+6\mathbf{e}^{-s_{1}}+6\mathbf{e}^{-s_{2}}\right)  }}{2\left(  \mathbf{e}%
^{-s_{1}}+\mathbf{e}^{-s_{2}}+\mathbf{e}^{s_{1}+s_{2}}\right)  }$\\
$z=\theta_{3}=-s_{1}-s_{2}+\log\frac{-\mathbf{e}^{s_{1}}-\mathbf{e}^{s_{2}%
}-\mathbf{e}^{-s_{1}-s_{2}}+\sqrt{\left(  \mathbf{e}^{2s_{1}}+\mathbf{e}%
^{2s_{2}}+6\mathbf{e}^{s_{1}+s_{2}}+\mathbf{e}^{-2(s_{1}+s_{2})}%
+6\mathbf{e}^{-s_{1}}+6\mathbf{e}^{-s_{2}}\right)  }}{2\left(  \mathbf{e}%
^{-s_{1}}+\mathbf{e}^{-s_{2}}+\mathbf{e}^{s_{1}+s_{2}}\right)  }$%
\end{tabular}
\ \ \right.  ,
\]
whose graphic representation is%
%TCIMACRO{\FRAME{dtbpFUX}{9.0325cm}{6.0231cm}{0pt}{\Qcb{figure 2 : The boundary
%of $D\left(  \mu_{\lambda}\right)  $}}{\Qlb{fig2}}{Plot}%
%{\special{ language "Scientific Word";  type "MAPLEPLOT";  width 9.0325cm;
%height 6.0231cm;  depth 0pt;  display "USEDEF";  plot_snapshots TRUE;
%mustRecompute FALSE;  lastEngine "MuPAD";  xmin "-10";  xmax "10";
%ymin "-10";  ymax "10";  xviewmin "-7";  xviewmax "3";  yviewmin "-7";
%yviewmax "3";  zviewmin "-7";  zviewmax "3";  viewset "XYZ";  rangeset "XYZ";
%phi 75;  theta -42;  cameraDistance "5.25759";
%cameraOrientation "[0,0,0.0419778]";  cameraOrientationFixed TRUE;
%plottype 5;  axesFont "Times New Roman,12,0000000000,useDefault,normal";
%num-x-gridlines 25;  num-y-gridlines 25;  plotstyle "wireframe";
%axesstyle "normal";  axestips FALSE;  plotshading "NONE";  lighting 0;
%xis \TEXUX{x};  yis \TEXUX{y};  var1name \TEXUX{$x$};  var2name \TEXUX{$y$};
%function \TEXUX{$G\left( u,v\right) $};  linestyle 1;  pointstyle "point";
%linethickness 1;  lineAttributes "Solid";  var1range "-10,10";
%var2range "-10,10";  surfaceColor "[flat::RGB:0000000000:0x000000ff]";
%surfaceStyle "Wire Frame";  num-x-gridlines 76;  num-y-gridlines 76;
%surfaceMesh "Mesh";  rangeset "XY";  VCamFile 'QTXU8400.xvz';valid_file "T";
%tempfilename 'QTXU8401.wmf';tempfile-properties "XPR";}} }%
%BeginExpansion

\begin{center}
\includegraphics{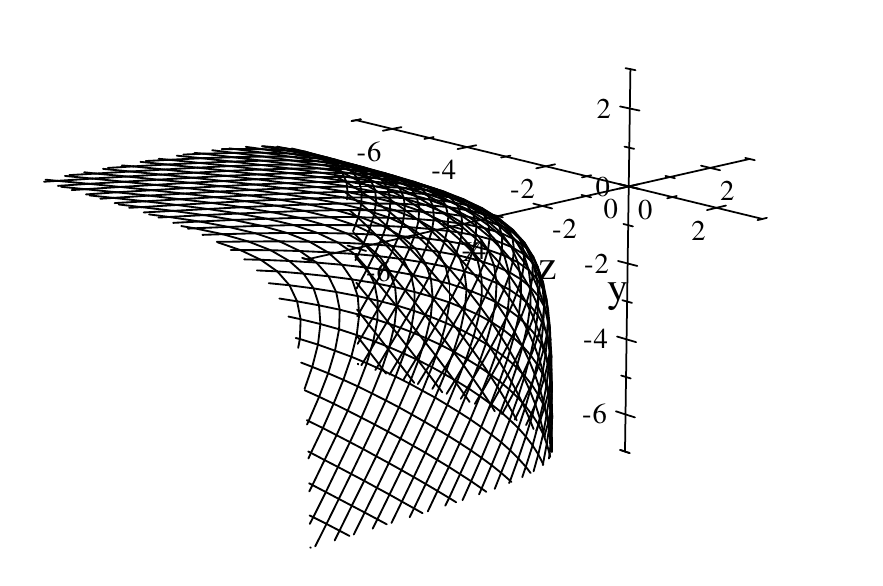}
\\
figure 2 : The boundary of $D\left(  \mu_{\lambda}\right)  $%
\label{fig2}%

\end{center}

%EndExpansion

If we choose $s_{1}=s_{2}=s_{3}=0,$ the condition $\overline{\theta}_{3}<\log
R_{\mathbf{s}}$ becomes $\overline{\theta}_{3}<\log\left(  \frac{1}{6}%
\sqrt{21}-\frac{1}{2}\right)  .$ As $\overline{\theta}_{3}=-\log4<\log\left(
\frac{1}{6}\sqrt{21}-\frac{1}{2}\right)  $, then $\left(  -\log4,-\log
4,-\log4\right)  \in D\left(  \mu_{\lambda}\right)  $.

Then, the introduction proves that
\[
\left(  \frac{A\left(  \frac{1}{4}z_{1},\frac{1}{4}z_{2},\frac{1}{4}%
z_{3}\right)  }{A\left(  \frac{1}{4},\frac{1}{4},\frac{1}{4}\right)  }\right)
^{-\lambda}=\left(  16-4z_{1}-4z_{2}-4z_{3}-z_{1}z_{2}-z_{1}z_{3}-z_{2}%
z_{3}\right)  ^{-\lambda}%
\]
is a generating function for all $\lambda>0$.

\begin{remark}
In this case, using the Remark \ref{mainremark} we obtain that an other
definition of $\mathbf{\theta}=\left(  \theta_{1},\theta_{2},\theta
_{3}\right)  $ $\in$ the boundary of $D\left(  u_{\lambda}\right)  $ is :
\[
\theta_{1}<0,\theta_{2}<-\log\left(  \frac{1-\mathbf{e}^{\theta_{1}}%
}{1+\mathbf{e}^{\theta_{1}}}\right)  ,\theta_{3}=\log\left(  \frac
{1-\mathbf{e}^{\theta_{1}}-\mathbf{e}^{\theta_{2}}-\mathbf{e}^{\theta
_{1}+\theta_{2}}}{1+\mathbf{e}^{\theta_{1}}+\mathbf{e}^{\theta_{2}}}\right)
,
\]
and an other parametrization of the boundary of $D\left(  u_{\lambda}\right)
$ is
\[
\left\{
\begin{tabular}
[c]{l}%
$x=\theta_{1}=u,u<0,$\\
$y=\theta_{2}=v-\log\left(  \frac{1-\mathbf{e}^{\theta_{1}}}{1+\mathbf{e}%
^{\theta_{1}}}\right)  ,v<0$\\
$z=\theta_{3}=\log\left(  \frac{1-\mathbf{e}^{\theta_{1}}-\mathbf{e}%
^{\theta_{2}}-\mathbf{e}^{\theta_{1}+\theta_{2}}}{1+\mathbf{e}^{\theta_{1}%
}+\mathbf{e}^{\theta_{2}}}\right)  $%
\end{tabular}
\ \right.  ,
\]
i.e.%
\[
\left\{
\begin{tabular}
[c]{l}%
$x=\theta_{1}=u,u<0,$\\
$y=\theta_{2}=v-\log\left(  \frac{1-\mathbf{e}^{u}}{1+\mathbf{e}^{u}}\right)
,v<0$\\
$z=\theta_{3}=\log\left(  \frac{\left(  1-\mathbf{e}^{2u}\right)  \left(
1-\mathbf{e}^{v}\right)  }{1+2\mathbf{e}^{u}+\mathbf{e}^{v}-\mathbf{e}%
^{u}\mathbf{e}^{v}+\mathbf{e}^{2u}}\right)  $%
\end{tabular}
\ \right.  .
\]
In addition, because $D\left(  \mu_{\lambda}\right)  $ is a convex set,
$\mathbf{\theta}=\left(  \theta_{1},\theta_{2},\theta_{3}\right)  \in D\left(
\mu_{\lambda}\right)  $ is defined by
\[
\theta_{1}<0,\theta_{2}<-\log\left(  \frac{1+\mathbf{e}^{\theta_{1}}%
}{1-\mathbf{e}^{\theta_{1}}}\right)  ,\theta_{3}<\log\left(  \frac
{1-\mathbf{e}^{\theta_{1}}-\mathbf{e}^{\theta_{2}}-\mathbf{e}^{\theta
_{1}+\theta_{2}}}{1+\mathbf{e}^{\theta_{1}}+\mathbf{e}^{\theta_{2}}}\right)
.
\]
Hence $\left(  -\log4,-\log4,-\log4\right)  \in D\left(  \mu_{\lambda}\right)
$ because $-\log4<0,-\log4<-\log\left(  \frac{1+\mathbf{e}^{-\log4}%
}{1-\mathbf{e}^{-\log4}}\right)  =-\ln\frac{5}{3}$ and $-\log4<$ $\log\left(
\frac{1-\mathbf{e}^{-\log4}-\mathbf{e}^{-\log4}-\mathbf{e}^{-\log4-\log4}%
}{1+\mathbf{e}^{-\log4}+\mathbf{e}^{-\log4}}\right)  =-\log\left(
4-4/7\right)  .$
\end{remark}

\section{Acknowledgment}

I thank G\'{e}rard Letac and Evelyne Bernadac for many helpful conversations.

\end{document}